\documentclass{article}
\usepackage{amssymb, amsmath}
\usepackage{microtype}
\usepackage{graphics}
\input xy
\xyoption{all}

\begin{document}
\title{Symplectic Killing spinors}
\author{Svatopluk Kr\'ysl \footnote{{\it E-mail address}: krysl@karlin.mff.cuni.cz}\\ {\it \small  Charles University, Sokolovsk\'a 83, Praha 8, Czech Republic.}
\thanks{The author of this article was supported by the grant GA\v{C}R 306-33/80397 of the Grant Agency of the Czech Republic. The work is a part of the research project MSM0021620839 financed by M\v{S}MT \v{C}R.}}
\maketitle

\noindent\centerline{\large\bf Abstract}  
 
 Let $(M,\omega)$ be a symplectic manifold admitting a metaplectic structure (a symplectic analogue of the Riemannian spin structure) and a torsion-free symplectic connection $\nabla.$ Symplectic Killing spinor fields for this structure are  sections of the symplectic spinor bundle satisfying a certain first order partial differential equation and they are the main topic of this paper. 
   We derive a necessary condition satisfied by a symplectic Killing spinor field. The advantage of this condition consists in the fact that it is expressed by a zeroth order operator.  This condition helps us substantionally to compute the symplectic Killing spinor fields for the standard symplectic vector spaces and the round sphere $S^2$ equipped with the volume form of the round metric.
 
{\it Math. Subj. Class.:} 58J60, 53C07



\section{Introduction}

In this article we shall study the so called symplectic Killing spinor fields on Fedosov manifolds admitting a metaplectic structure.  A Fedosov manifold 
is a  structure  consisting of a symplectic manifold $(M^{2l},\omega)$ and the so called Fedosov connection on $(M,\omega)$. A Fedosov connection $\nabla$ is an affine connection on $(M,\omega)$ such that it is
 symplectic, i.e., $\nabla \omega =0,$ and torsion-free. Let us notice that in contrary to the Riemannian geometry, a  Fedosov connection is not unique. Thus, it seems natural to add the Fedosov connection into the studied structure obtaining the notion of a Fedosov manifold in this way. See, e.g., Tondeur \cite{Tondeur} for symplectic connections for presymplectic structures and Gelfand, Retakh, Shubin \cite{GRS} for Fedosov connections.  

It is known that if $l>1,$ the curvature tensor of a Fedosov connection decomposes into two invariant parts, namely into the so called symplectic Ricci curvature and symplectic Weyl curvature tensor fields. If $l=1,$ only the symplectic Ricci curvature occurs. See Vaisman \cite{Vaisman} for details.

 In order to define a symplectic Killing spinor field, we shall briefly describe the so called metaplectic structures with help of which these fields are defined.
Any symplectic group $Sp(2l,\mathbb{R})$ admits a non-trivial, i.e., connected, two-fold covering, the so called metaplectic group, denoted by
$Mp(2l,\mathbb{R})$ in this paper. A mataplectic structure over a symplectic manifold is a symplectic analogue of the Riemannian spin structure. 
In particular, one of its parts is a principal $Mp(2l,\mathbb{R})$-bundle which covers twice the bundle of symplectic reperes of $(M^{2l},\omega).$  Let us denote this principal $Mp(2l,\mathbb{R})$-bundle by $q:~\mathcal{Q} \to M.$

 Now, let us say a few words about  the symplectic spinor fields.
These  fields are sections of the so called symplectic spinor bundle $\mathcal{S} \to M.$ This vector bundle is the  bundle associated to the principal $Mp(2l,\mathbb{R})$-bundle $q:\mathcal{Q}\to M$ via the so called {\it Segal-Shale-Weil representation}. The Segal-Shale-Weil representation is a distinguished representation of the metaplectic group and plays a similar role in the quantization of boson particles as the  spinor representations of spin groups play in the quantization of  fermions. See, e.g., Shale \cite{Shale}. The Segal-Shale-Weil representation is unitary and does not descend to a representation of the symplectic group.  The vector space of the underlying Harish-Chandra $(\mathfrak{g},K)$-module of the Segal-Shale-Weil representation is isomorphic to  $S^{\bullet}(\mathbb{R}^l),$ the {\it symmetric} power of a Lagrangian subspace $\mathbb{R}^l$ in the symplectic vector space $\mathbb{R}^{2l}.$ Thus, the situation is parallel to the complex orthogonal case, where the spinor representation can be realized on the  {\it exterior} power of a maximal isotropic subspace.
The Segal-Shale-Weil representation and some of its analytic versions are sometimes called oscillatory representation, metaplectic representation or symplectic spinor representation. For a more detailed explanation of 
the last name, see, e.g., Kostant \cite{Kostant}. 

 The symplectic Killing spinor field is a non-zero section of the symplectic spinor bundle $\mathcal{S}\to M$ satisfying certain linear first order partial differential equation formulated by the connection $\nabla^S:\Gamma(M,\mathcal{S})\times \Gamma(M,TM) \to \Gamma(M,\mathcal{S}),$ the associated connection to the Fedosov connection $\nabla.$ 
This partial differential equation is a symplectic analogue of the classical symplectic Killing spinor equation from at least two aspects.
One of them is rather formal. Namely, the defining equation for a symplectic Killing spinor is of the "same shape" as that one for a Killing spinor field on a Riemannian spin manifold. The second similarity can be expressed by comparing this
equation with the so called symplectic Dirac equation and the symplectic twistor equation and will be discussed bellow 
in this paper. Let us mention that   any symplectic Killing spinor field determines a unique complex number, the so called symplectic Killing spinor number. 
Notice that the symplectic Killing spinor fields were considered already in a connection with the existence of a linear embedding of the spectrum
of the so called symplectic Dirac operator into the spectrum of the so called symplectic Rarita-Schwinger operator. The symplectic Killing spinor fields represent an obstruction for the mentioned embedding. See Kr\'ysl \cite{KryslRarita} for this aspect.

In many particular cases, the equation for symplectic Killing spinor fields seems  to be rather complicated. On the other hand, in many cases it is known that its solutions are rare. Therefore it is reasonable to look for a  necessary condition satisfied by a symplectic Killing spinor field which is simpler than the defining equation itself. 
Let us notice that similar necessary conditions are known and parallel methods were used in   Riemannian or Lorentzian spin geometry. See, e.g., Friedrich \cite{Friedrich}. 

In this paper, we shall prove that any symplectic Killing spinor field necessarilly satisfies  certain zeroth order differential equation. More precisely, we prove that any symplectic Killing spinor is necessarily a section of the kernel of a symplectic spinor bundle morphism.  We derive  this equation by prolongating the symplectic Killing spinor equation.  We make such a prolongation  that enables us to compare the result with an appropriate  part of the curvature tensor of the associated connection $\nabla^S$ acting on symplectic spinors. An explicit formula for this part of the curvature action was already derived in Kr\'ysl \cite{Kryslcurv}. Especially, it is known that the symplectic Weyl curvature   of $\nabla$ does not show up in this part and thus, the mentioned morphism depends  on the symplectic Ricci part of the curvature of the Fedosov connection $\nabla$ only.
This will make us able to prove that the only symplectic Killing number of a Fedosov manifold of  Weyl type is zero.
This will in turn  imply that any symplectic Killing spinor on the standard symplectic vector space of an arbitrary finite dimension and equipped with the standard flat connection is constant.  This result can be obtained easily when one knows the prolongated equation, whereas computing the symplectic Killing spinors without this knowledge is rather complicated. This fact  will be illustrated when we will compute the symplectic Killing spinors on the standard symplectic $2$-plane using just the defining equation for symplectic Killing spinor field.
 
 The cases when the prolongated equation does not help so easily as in the case of the Weyl type Fedosov manifolds are the Ricci type ones.
Nevertheless, we prove that there are no symplectic Killing spinors on the $2$-sphere, equipped with the volume form of the  round metric as the symplectic form and the Riemannian connection as the Fedosov connection.
Let us remark that in this case, the prolongated equation  has a shape of a stationary Schr\"odinger equation. More precisely, it has the shape of
the equation for the eigenvalues of certain  oscillator-like quantum Hamiltonian determined  completely by the symplectic Ricci curvature tensor of the Fedosov connection.

Let us notice that there are some applications of symplectic spinors in physics besides those in the mentioned article of Shale \cite{Shale}. For an application in string theory physics, see, e.g., Green, Hull \cite{GH}.
  
 In the second section, some necessary notions from symplectic linear algebra and representation theory of reductive Lie groups are explained and the Segal-Shale-Weil 
 representation and the symplectic Clifford multiplication are introduced. In the third section, the Fedosov connections are introduced and some properties of their curvature tensors acting on symplectic spinor fields are summarized. In the fourth section, the symplectic Killing spinors are defined and symplectic Killing spinors on the standard symplectic $2$-plane are computed. In this section, a connection of the symplectic Killing spinor fields to the eigenfunctions of symplectic Dirac and symplectic twistor operators is formulated and proved. Further, the mentioned prolongation of the symplectic Killing spinor equation is derived  and the symplectic Killing spinor fields on the standard symplectic vector spaces are computed. At the end, the case of the round  sphere $S^2$  is treated.

\section{Symplectic spinors and symplectic spinor valued forms}

Let us start recalling some notions from  symplectic linear algebra. Let us mention that we shall often use the Einstein summation convention without mentioning it explicitly.
Let $(\mathbb{V},\omega_0)$ be a symplectic vector space of dimension $2l,$ i.e., $\omega_0$ is a non-degenerate anti-symmetric bilinear form on the vector space $\mathbb{V}$ of dimension $2l.$ Let $\mathbb{L}$ and $\mathbb{L}'$ be two Lagrangian subspaces\footnote{i.e., maximal isotropic wr. to $\omega_0,$ in particular $\dim \mathbb{L} = \dim \mathbb{L}'=l$} of $(\mathbb{V},\omega_0)$ such that $\mathbb{L} \oplus \mathbb{L}' =\mathbb{V}.$
Let $\{e_i\}_{i=1}^{2l}$ be an adapted symplectic basis of $(\mathbb{V}=\mathbb{L}\oplus \mathbb{L}',\omega_0),$ i.e.,
$\{e_i\}_{i=1}^{2l}$ is a symplectic basis and $\{e_i\}_{i=1}^l\subseteq \mathbb{L}$ and $\{e_i\}_{i=l+1}^{2l} \subseteq \mathbb{L}'.$ Because the definition of a symplectic basis is not unique, we shall fix one which we shall use in this text. A basis $\{e_i\}_{i=1}^{2l}$ of $(\mathbb{V},\omega_0)$ is called symplectic, if
$\omega_0(e_i, e_j)=1$ iff $1\leq i \leq l$ and $j=l+i;$ $\omega_0(e_i,e_j)=-1$ iff $l+1\leq i\leq 2l$ and $j=i-l$ and $\omega_0(e_i,e_j)=0$ in the remaining cases. Whenever a symplectic basis will be chosen, we will denote the basis of $\mathbb{V}^*$ dual to $\{e_i\}_{i=1}^{2l}$ by $\{\epsilon^i\}_{i=1}^{2l}.$
Further for $i,j=1,\ldots, 2l$, we set $\omega_{ij}:=\omega_0(e_i,e_j)$ an similarly for other type of tensors. For $i,j=1,\ldots, 2l,$ we define $\omega^{ij}$ by the equation $\sum_{k=1}^{2l}\omega_{ik}\omega^{jk}=\delta^i_j.$

  As in the orthogonal case, we would like to rise and lower indices.
Because the symplectic form $\omega_0$ is antisymmetric, we should be more careful in this case. 
For coordinates ${K_{ab\ldots c\ldots d}}^{rs \ldots t \ldots u}$ of a tensor $K$ over $\mathbb{V},$  we denote
the expression $\omega^{ic}{K_{ab\ldots c \ldots d}}^{rs \ldots t}$ by 
${{{K_{ab \ldots}}^{i}}_{\ldots d}}^{rs \ldots t}$ and 
${K_{ab\ldots c}}^{rs \ldots t \ldots u}\omega_{ti}$ by ${{{K_{ab \ldots c}}^{rs\ldots}}_{i}}^{\ldots u}$ and similarly for other types of tensors and also in a geometric setting when we will be considering tensor fields over a symplectic manifold $(M,\omega)$.   

Let us denote the symplectic group $Sp(\mathbb{V},\omega_0)$ of $(\mathbb{V},\omega_0)$ by $G.$
Because the maximal compact subgroup of $G$ is isomorphic to the unitary group $U(l)$ which is of homotopy type $\mathbb{Z},$ we have  $\pi_1(G)\simeq \mathbb{Z}.$ From the theory of covering spaces, we know that there exists up to an isomorphism a unique connected double cover of $G.$ This double cover is the so called metaplectic group $Mp(\mathbb{V},\omega_0)$ and will be denoted by $\tilde{G}$ in this text.  We shall denote the covering homomorphism by $\lambda,$ i.e.,
$\lambda:\tilde{G} \to G$ is a fixed member of the isomorphism class of all connected 2:1 coverings.

Now, let us  recall some notions from representation theory of reductive Lie groups which we shall need in this paper. Let us mention that these notions are rather of technical character for the purpose of our article.
For a reductive Lie group $G$ in the sense of Vogan \cite{Vogan}, let $\mathcal{R}(G)$ be the category the object of which are  complete, locally convex, Hausdorff vector spaces with a continuous action of $G$ which is admissible and of finite length; the morphisms are continuous linear $G$-equivariant maps between the objects. Let us notice that, e.g., finite covers of the classical groups are reductive. It is known  that any irreducible unitary representation of a reductive group is in $\mathcal{R}(G).$ Let $\mathfrak{g}$ be the Lie algebra of $G$ and $K$ be a maximal compact subgroup of $G.$ It is well known that there exists the so called $L^2$-globalization functor, denoted by $L^2$ here, from the category of admissible Harish-Chandra modules to the category $\mathcal{R}(G)$. See Vogan \cite{Vogan} for details. Let us notice that this functor behaves compatibly wr. to Hilbert tensor products. See, e.g., Vogan \cite{Vogan} again.  For an object $\bf E$ in $\mathcal{R}(G),$ let us denote its underlying Harish-Chandra $(\mathfrak{g},K)$-module by $E$ and when we will be considering only its $\mathfrak{g}^{\mathbb{C}}$-module structure, we shall denote it by $\mathbb{E}.$ 
If $\mathfrak{g}^{\mathbb{C}}$ happens to be a simple complex Lie algebra of rank $l$, let us denote its Cartan subalgebra  by $\mathfrak{h}^{\mathbb{C}}.$ The set $\Phi$ of roots for $(\mathfrak{g}^{\mathbb{C}}, \mathfrak{h}^{\mathbb{C}})$ is then uniquely determined.  Further let us choose a set $\Phi^+\subseteq \Phi$ of positive roots and denote the corresponding set of fundamental weights by $\{\varpi_i\}_{i=1}^l.$ 
For $\lambda \in \mathfrak{h}^{\mathbb{C}},$ let us denote the irreducible highest weight module
with the highest weight $\lambda$ by $L(\lambda).$
 
Let us denote by $\mathcal{U}({\bf W})$ the group of unitary operators on a Hilbert space $\bf W$ and let
$L: Mp(\mathbb{V},\omega_0) \to \mathcal{U}({\bf L}^2(\mathbb{L}))$ be the Segal-Shale-Weil representation of the 
metaplectic group. 
 It is an infinite dimensional unitary representation of the metaplectic group on the complex valued square Lebesgue integrable functions defined on the Lagrangian subspace $\mathbb{L}.$ This representation does  not descend to a representation of the symplectic group $Sp(\mathbb{V},\omega_0).$ See, e.g., Weil \cite{Weil} and Kashiwara, Vergne \cite{KV}. 
For our convenience, let us set ${\bf S}:={\bf L}^2(\mathbb{L})$ and call it the symplectic spinor module and its elements symplectic spinors.
It is well known that as a $\tilde{G}$-module, $\bf S$ decomposes into the direct sum ${\bf S} ={\bf S}_+ \oplus {\bf S}_-$ of two irreducible submodules. The submodule ${\bf S}_+$ (${\bf S}_-$) consists of even (odd) functions in ${\bf L}^2(\mathbb{L}).$
 Further, let us notice that in K\'ysl \cite{KryslJRT}, a slightly different analytic version (based on the so called minimal globalizations) of this representation was used.


As in the orthogonal case, we may multiply spinors by vectors. The multiplication $.: \mathbb{V} \times {\bf S} \to {\bf S}$ will be called symplectic Clifford multiplication and it is defined as follows.
For $f \in {\bf S}$ and $i=1,\ldots,l,$ we set
$$(e_i.f)(x):=\imath x^if(x),$$
$$(e_{l+i}.f)(x):=\frac{\partial f}{\partial x^i}(x), \, x \in \mathbb{L}$$
and extend it linearly to get the symplectic Clifford multiplication.
The  symplectic Clifford multiplication (by a fixed vector) has to be understood as an unbounded operator on ${\bf L}^2(\mathbb{L}).$ See Habermann, Habermann \cite{HH} for details. 
Let us also notice that the symplectic Clifford multiplication corresponds to the so called Heisenberg canonical quantization known from quantum mechanics. (For brevity, we shall write $v.w.s,$ instead of $v.(w.s)$, $v,w \in \mathbb{V}$ and $s\in {\bf S}.$)

It is easy to check that the symplectic Clifford multiplication satisfies the relation described in the following

{\bf Lemma 1:} For $v, w\in \mathbb{V}$ and $s \in {\bf S},$
 we have  $$v.(w.s)-w.(v.s)=-\imath \omega_0(v,w)s.$$

{\it Proof.} See Habermann, Habermann \cite{HH}. $\Box$

 Let us consider the representation $$\rho: \tilde{G} \to \mbox{Aut}(\bigwedge^{\bullet}\mathbb{V}^*\otimes {\bf S})$$ of the metaplectic group $\tilde{G}$ on $\bigwedge^{\bullet}\mathbb{V}^*\otimes {\bf S}$ given by 
$$\rho(g)(\alpha \otimes s):=\lambda^{* \wedge r}(g)\alpha \otimes L(g)s,$$ where $r=0,\ldots, 2l,$ $\alpha \in  \bigwedge^r\mathbb{V}^*,$ $s \in {\bf S}$ and $\lambda^{*\wedge r}$  denotes the $r$th wedge power of the representation $\lambda^*$ dual to $\lambda,$ and extended linearly.  For definiteness, let us consider the vector space $\bigwedge^{\bullet}\mathbb{V}^*\otimes {\bf S}$   equipped with the topology of the Hilbert tensor product. Because the $L^2$-globalization functor behaves compatibly wr. to the Hilbert tensor products, one can easily see that the representation $\rho$ belongs to the class $\mathcal{R}(\tilde{G}).$ 

In the next theorem, the space o symplectic valued exterior  two-forms is decomposed into irreducible summands.

{\bf Theorem 2:} For $\frac{1}{2}\mbox{dim}(\mathbb{V}) =:l > 2$, the following isomorphisms
$$\bigwedge^2\mathbb{V}^* \otimes {\bf S}_{\pm} \simeq {\bf E}^{20}_{\pm} \oplus  {\bf E}^{21}_{\pm}
\oplus {\bf E}^{22}_{\pm}$$ hold. 
For   $j_2=0,1,2,$ the
 ${\bf E}^{2j_2}$ are uniquely determined by the conditions that first, they are submodules of the corresponding tensor products and second, 
\begin{eqnarray*}
\mathbb{E}^{20}_{-} \simeq \mathbb{S}_{-} \simeq L(\varpi_{l-1}-\frac{3}{2}\varpi_l), \, \mbox{  } \mathbb{E}^{20}_+ \simeq \mathbb{S}_+ \simeq L(-\frac{1}{2}\varpi_{l}),\\
\mathbb{E}^{21}_- \simeq  L(\varpi_1-\frac{1}{2}\varpi_l), \, \mbox{   } \mathbb{E}^{21}_+ \simeq L(\varpi_1+\varpi_{l-1}-\frac{3}{2}\varpi_l),\\ \mathbb{E}^{22}_+  \simeq  L(\varpi_{2}-\frac{1}{2}\varpi_l) \, \mbox{ and } \, \mathbb{E}^{22}_- \simeq  L(\varpi_2+\varpi_{l-1}-\frac{3}{2}\varpi_l).
\end{eqnarray*}
{\it Proof.} This theorem is proved in Kr\'ysl \cite{KryslRarita} or Kr\'ysl \cite{KryslJRT} for the so called minimal globalizations. Because the $L^2$-globalization behaves compatibly wr. to the considered Hilbert tensor product topology, the statement remains true. $\Box$

{\bf Remark:} Let us notice that for $l=2,$ the number of irreducible summands in symplectic spinor valued two-forms is the same as that one for $l>2.$
  In this case ($l=2$), one only has   to change the prescription for the highest weights described in the preceding theorem. For $l=1,$ the number of the irreducible summands is different  from that one for $l\geq 2.$ Nevertheless, in this case the decomposition is also multiplicity-free. See Kr\'ysl \cite{KryslJRT} for details.

In order to make some proofs in the section on symplectic Killing spinor fields simpler and more clear, let us introduce the operators

\begin{eqnarray*}
F^+&:& \bigwedge^{\bullet}\mathbb{V}^*\otimes {\bf S} \to
\bigwedge ^{\bullet+1}\mathbb{V}^*\otimes {\bf S}, \, \mbox{  }F^+(\alpha \otimes
s):=\sum_{i=1}^{2l}\epsilon^i\wedge \alpha \otimes e_i.s,\\
F^-&:& \bigwedge^{\bullet}\mathbb{V}^{*}\otimes {\bf S} \to \bigwedge^{\bullet-1}\mathbb{V}^*\otimes {\bf S}, \, \mbox{  } F^-(\alpha \otimes s):=-\sum_{i,j=1}^{2l} \omega^{ij}\iota_{e_i}\alpha \otimes e_j.s,\\
H&:&\bigwedge^{\bullet}\mathbb{V}^*\otimes {\bf S} \to \bigwedge^{\bullet}\mathbb{V}^*\otimes {\bf S}, \, \mbox{  } H:=\{F^+,F^-\}.
\end{eqnarray*}

{\bf Remark:} 
\begin{itemize}
\item[1)] One easily finds out that the operators are independent of the choice of an adapted symplectic basis $\{e_i\}_{i=1}^{2l}.$
\item[2)] Let us remark that the operators $F^+, F^-$ and $H$ defined here differ from the operators $F^+, F^-, H$ defined in Kr\'ysl \cite{KryslJRT}. Though, by a constant real multiple only.
\item[3)]
The operators $F^+$ and $F^-$ are used to prove the Howe correspondence for $Mp(\mathbb{V},\omega_0)$ acting on 
$\bigwedge^{\bullet}\mathbb{V}^*\otimes {\bf S}$ via the representation $\rho.$ More or less, the ortho-symplectic super Lie algebra $\mathfrak{osp}(1|2)$ plays the role of a (super Lie) algebra, a representation of which is the appropriate commutant. See Kr\'ysl \cite{KryslJRT} for details. 
\end{itemize}

In the next lemma the $\tilde{G}$-equivariance of the operators $F^+, F^-$ and $H$ is stated, some properties of
$F^{\pm}$ are mentioned and the value of 
$H$ on degree-homogeneous elements is computed. We shall use this lemma when we will be  treating the symplectic Killing spinor fields in the fourth section.

{\bf Lemma 3:} Let $(\mathbb{V}=\mathbb{L}\oplus \mathbb{L}',\omega_0)$ be a $2l$ dimensional
symplectic vector space. Then
\begin{itemize}
\item[1)] the operators $F^+,$ $F^+$ and $H$ are $\tilde{G}$-equivariant,
\item[2)] 
\begin{itemize}
 \item[a)]$F^-_{|\bf{E}^{11}}=0,$ 
 \item[b)]$F^+_{|{\bf E}^{00}}$ is an isomorphism onto ${\bf E}^{10}$   
 \item[c)]$(F^+)^2_{|{\bf S}}=-\frac{\imath}{2}\omega \otimes \mbox{Id}_{|{\bf S}}$ and it is an isomorphism onto ${\bf E}^{20}.$ 
\end{itemize}
\item[3)] For  $r=0,\ldots, 2l,$ we have
$$H_{|\bigwedge^r\mathbb{V}^* \otimes {\bf S}}=\imath(r-l)\hbox{Id}_{|\bigwedge^r\mathbb{V}^*\otimes {\bf S}}.$$
\end{itemize}
{\it Proof.} See Kr\'ysl \cite{KryslJRT}. $\Box$

Let us remark that the items 1 and 3 of the preceding lemma follow by a direct computation, and the second item follows from the first item, decomposition theorem (Theorem 2), a  version of the Schur lemma and a direct computation.
 
\section{Curvature of Fedosov manifolds and its actions on symplectic spinors}

After we have finished the "algebraic part" of this paper, let us recall some basic facts on Fedosov manifolds, their curvature tensors, metaplectic structures and the action of the curvature tensor on symplectic spinor fields.

Let us start recalling some notions and results related to the so called  Fedosov manifolds.
Let $(M^{2l},\omega)$ be a symplectic manifold of dimension $2l.$ Any torsion-free affine connection $\nabla$ on $M$ preserving $\omega,$ i.e., $\nabla \omega =0,$  is called {\it Fedosov connection}. Let us recall that torsion-free means $T(X,Y):=\nabla_XY-\nabla_YX -[X,Y]=0$ for all vector fields  $X,Y \in \mathfrak{X}(M).$ 
The triple $(M,\omega,\nabla),$ where $\nabla$ is a Fedosov connection, will be called {\it Fedosov manifold}.
As we have already mentioned in the Introduction,  a Fedosov connection for a given symplectic manifold $(M,\omega)$ is not unique. 
Let us remark that Fedosov manifolds are used for a construction of geometric quantization of symplectic manifolds due to Fedosov. See, e.g., Fedosov \cite{Fedosov}.

To fix our notation, let us recall the classical definition of the curvature tensor $R^{\nabla}$ of the connection $\nabla,$ we shall be using here. We set 
$$R^{\nabla}(X,Y)Z:=\nabla_X\nabla_Y Z - \nabla_Y \nabla_X Z - \nabla_{[X,Y]}Z$$ for
$X,Y,Z \in \mathfrak{X}(M).$ 

Let us choose a local adapted symplectic frame $\{e_i\}_{i=1}^{2l}$ on a fixed open subset $U\subseteq M.$
By a local adapted  symplectic frame $\{e_i\}_{i=1}^{2l}$ over $U,$ we mean such a local frame that
for each $m \in U$ the basis $\{(e_i)_m\}_{i=1}^{2l}$ is an adapted symplectic basis of $(T_mM,\omega_m).$
Whenever a symplectic frame is chosen, we denote its dual frame   by $\{\epsilon^i\}_{i=1}^{2l}.$
Although some of the formulas bellow hold only locally, containing a local adapted symplectic frame, we will not mention this restriction explicitly.
  
From the symplectic curvature tensor field $R^{\nabla}$, we can build the symplectic Ricci curvature tensor field  $\sigma^{\nabla}$ defined by the 
classical formula
$$\sigma^{\nabla}(X,Y):=\mbox{Tr}(V \mapsto R^{\nabla}(V,X)Y)$$ for each $X,Y \in \mathfrak{X}(M)$ (the variable $V$ denotes a vector field on $M$). For the chosen frame and $i,j=1,\ldots, 2l$, we define
$$\sigma_{ij}:=\sigma^{\nabla}(e_i,e_j).$$


Let us define the extended Ricci tensor field by the equation
$$\widetilde{\sigma}(X,Y,Z,U):=\widetilde{\sigma}_{ijkn}X^iY^jZ^kU^n, \mbox{ } X,Y,Z,U \in \mathfrak{X}(M),$$
where for $i,j,k,n=1,\ldots, 2l,$
$$2(l+1)\widetilde{\sigma}_{ijkn}:=\omega_{in}\sigma_{jk}-\omega_{ik}\sigma_{jn}+\omega_{jn}\sigma_{ik}-\omega_{jk}\sigma_{in}+2\sigma_{ij}\omega_{kn}.$$

A Fedosov manifold $(M,\omega,\nabla)$ is called of Weyl type, if $\sigma=0.$ Let us notice, that it is called of  
Ricci type, if $R=\widetilde{\sigma}.$ In Vaisman \cite{Vaisman}, one can find  more information on the $Sp(2l,\mathbb{R})$-invariant decomposition of the curvature tensors of  Fedosov connections.

Now, let us describe the geometric structure with help of which the  
symplectic  Killing  spinor fields  are defined.  This structure, called {\it metaplectic}, is a 
symplectic analogue of the notion of a spin structure in  the Riemannian geometry.
For a symplectic manifold $(M^{2l}, \omega)$  of dimension $2l,$
let us denote the bundle of symplectic reperes in $TM$ by
$\mathcal{P}$ and  the foot-point projection of $\mathcal{P}$ onto
$M$ by $p.$ Thus $(p:\mathcal{P}\to M, G),$ where $G\simeq
Sp(2l,\mathbb{R}),$ is a principal $G$-bundle over $M$. As in
the subsection 2, let $\lambda: \tilde{G}\to G$ be a member
of the isomorphism class of the non-trivial two-fold coverings of
the symplectic group $G.$ In particular, $\tilde{G}\simeq
Mp(2l,\mathbb{R}).$ Further, let us consider a principal
$\tilde{G}$-bundle $(q:\mathcal{Q}\to M, \tilde{G})$ over the
symplectic manifold $(M,\omega).$ We call a pair
$(\mathcal{Q},\Lambda)$   metaplectic structure if  $\Lambda:
\mathcal{Q} \to \mathcal{P}$ is a surjective bundle homomorphism
over the identity on $M$ and if the following diagram,
$$\begin{xy}\xymatrix{
\mathcal{Q} \times \tilde{G} \ar[dd]^{\Lambda\times \lambda} \ar[r]&   \mathcal{Q} \ar[dd]^{\Lambda} \ar[dr]^{q} &\\
                                                            & &M\\
\mathcal{P} \times G \ar[r]   & \mathcal{P} \ar[ur]_{p} }\end{xy}$$
with the
horizontal arrows being respective actions of the displayed groups, commutes.
See, e.g.,  Habermann,  Habermann \cite{HH} and Kostant \cite{Kostant} for
details on   metaplectic structures. Let us only remark that typical examples of symplectic manifolds admitting a metaplectic structure are cotangent bundles of orientable manifolds (phase spaces), Calabi-Yau manifolds and complex projective spaces $\mathbb{CP}^{2k+1}$, $k \in \mathbb{N}_0.$

Let us denote  the  vector bundle associated to the introduced principal $\tilde{G}$-bundle
$(q:\mathcal{Q}\to M,\tilde{G})$  via the representation $\rho$ acting on ${\bf S}$ by $\mathcal{S},$ and call this associated vector bundle {\it symplectic spinor bundle}.  Thus, we have $\mathcal{S}=\mathcal{Q}\times_{\rho}{\bf S}.$ The sections $\phi \in \Gamma(M,\mathcal{S})$ will be called {\it symplectic spinor fields}.
 Further for  $j_2=0,1, 2$, we define the  associated vector bundles $\mathcal{E}^{2j_2}$ by the prescription  $\mathcal{E}^{2j_2}:=\mathcal{Q}\times_{\rho} {\bf E}^{2j_2}.$   Further, we define
 $\mathcal{E}^r:=\Gamma(M,\mathcal{Q}\times_{\rho}\bigwedge^r\mathbb{V}^*\otimes {\bf S}),$ i.e., the space o symplectic spinor valued differential $r$-forms, $r=0,\ldots, 2l.$
Because the symplectic Clifford multiplication is $\tilde{G}$-equivariant (see Habermann, Habermann \cite{HH}), we can lift it to the associated vector bundle structure, i.e., to let it act on the elements from $\Gamma(M,\mathcal{S}).$ For  $j_2=0, 1, 2,$ let us denote the vector bundle projections $\Gamma(M,\mathcal{E}^2)\to \Gamma(M,\mathcal{E}^{2{j_2}})$ by $p_{2j_2},$ i.e.,
$p_{2j_2}:\Gamma(M,\mathcal{E}^2)\to \Gamma(M,\mathcal{E}^{2j_2})$ for all appropriate $j_2.$ This definition makes sense because due to the decomposition result (Theorem 2) and the Remark bellow the Theorem 2, the $\tilde{G}$-module of symplectic spinor valued exterior $2$-forms is multiplicity-free.

Let $Z$ be the principal bundle connection on the principal $G$-bundle $(p:\mathcal{P}\to M, G)$ associated to the chosen Fedosov connection $\nabla$ and $\tilde{Z}$ be a lift of $Z$ to the principal $\tilde{G}$-bundle $(q:\mathcal{Q}\to M,\tilde{G}).$ Let us denote by $\nabla^S$ the covariant derivative associated to $\tilde{Z}$. Thus, in particular, $\nabla^S$ acts on the symplectic spinor fields.

Any section $\phi$ of the associated vector bundle $\mathcal{S}=\mathcal{Q}\times_{\rho}{\bf S}$ can be equivalently considered as
a $\tilde{G}$-equivariant $\bf S$-valued function on $\mathcal{Q}.$ Let us denote this function by $\hat{\phi},$ i.e.,
$\hat{\phi}:\mathcal{Q} \to {\bf S}.$ For a local adapted symplectic frame $s:U \to \mathcal{P},$ let us denote by $\overline{s}:U \to \mathcal{Q}$ one of the lifts of $s$ to $\mathcal{Q}.$  Finally, let us set $\phi_s:=\hat{\phi}\circ \overline{s}.$
Further for $q \in \mathcal{Q}$ and $\psi \in {\bf S}$, let us denote by $[q,\psi]$ the equivalence class in $\mathcal{S}$ containing $(q,\psi).$ (As it is well known, the total space $\mathcal{S}$ of the  symplectic spinor bundle is the product  $\mathcal{Q} \times {\bf S}$ modulo an equivalence relation.)

{\bf Lemma 4:} Let $(M,\omega, \nabla)$ be a Fedosov manifold admitting a metaplectic structure.
Then for each $X \in \mathfrak{X}(M),$ $\phi \in \Gamma(M,\mathcal{S})$ and a local adapted symplectic frame $s: U \to \mathcal{P},$
we have
$$\nabla_X^S\phi = [\overline{s},X(\phi_s)]-\frac{\imath}{2}\sum_{i=1}^l[e_{i+l}.(\nabla_X e_i).-e_i.(\nabla_X e_{i+l}).]\phi \mbox{ and }$$
$$\nabla_X^S(Y.\phi)=(\nabla^S_XY) .\phi + X.\nabla_Y^S\phi.$$

{\it Proof.} See Habermann, Habermann \cite{HH}. $\Box$

The curvature tensor on symplectic spinor fields is defined by the formula 
$$R^S(X,Y)\phi=\nabla^S_X\nabla_Y^S\phi-\nabla_Y^S\nabla^S_X\phi-\nabla_{[X,Y]}^S\phi,$$ 
where $\phi \in \Gamma(M,\mathcal{S})$ and $X,Y \in \mathfrak{X}(M).$ 
In the next lemma, a part of the action of $R^S$ on the space of symplectic spinors is described using just the symplectic Ricci curvature tensor field $\sigma.$
 
{\bf Lemma 5:} Let $(M,\omega,\nabla)$ be a Fedosov manifold  admitting a metaplectic structure.
Then for a symplectic spinor field $\phi \in \Gamma(M,\mathcal{S}),$ we have
$$p^{20}R^S\phi=\frac{\imath}{2l}\sigma^{ij} \omega_{kl} \epsilon^k\wedge \epsilon^l \otimes e_i.e_j.\phi.$$

  
{\it Proof.} See  Kr\'ysl \cite{Kryslcurv}. $\Box.$

\section{Symplectic Killing spinor fields}

In this section, we shall focus our attention to the symplectic Killing spinor fields. More specifically, we compute the symplectic Killing spinor fields on some Fedosov manifolds admitting a metaplectic structure and derive a necessary condition satisfied by a symplectic Killing spinor field.

Let $(M,\omega,\nabla)$ be a Fedosov manifold admitting a metaplectic structure. We call a non-zero section $\phi \in \Gamma(M,\mathcal{S})$ {\it symplectic Killing spinor field} if
                 $$\nabla_X^S \phi = \lambda X. \phi$$ 
for a complex number $\lambda \in \mathbb{C}$ and each vector field $X \in \mathfrak{X}(M).$ 
The complex number $\lambda$ will sometimes be called symplectic Killing spinor number.
(Allowing the zero section to be a symplectic Killing spinor would make the notion of a symplectic Killing spinor number meaningless.)

Let us note that one can rewrite equivalently the preceding defining equation for a symplectic Killing spinor into  
$$\nabla^S\phi = \lambda F^+ \phi.$$ Indeed, if this equation is satisfied, we get by inserting
the local vector field $X=X^ie_i$ the equation $\nabla^S_X \phi = \iota_X (\lambda \epsilon^i \otimes e_i.\phi) = \lambda \epsilon^i(X) e_i.\phi = \lambda X^ie_i.\phi =\lambda X.\phi,$ i.e., the defining equation. Conversely, one can prove that $\nabla_X^S \phi =\lambda X.\phi$ is equivalent to $\iota_X \nabla^S \phi = \iota_X (\lambda F^+ \phi). $
Because this equation holds for each vector field $X,$ we get $\nabla^S \phi = \lambda F^+ \phi.$
We shall call  both the defining equation and the equivalent equation for a symplectic Killing spinor field the
symplectic Killing spinor equation.

In the next example, we compute the symplectic Killing spinors on the standard symplectic $2$-plane.

{\bf Example 1:} Let us solve the symplectic Killing spinor equation for the standard symplectic vector space $(\mathbb{R}^{2}[s,t],\omega_0)$ equipped with the standard flat Euclidean connection $\nabla.$ In this case, $(\mathbb{R}^{2},\omega_0,\nabla)$ is also a Fedosov manifold.
The bundle of symplectic reperes in $T\mathbb{R}^2$ defines a principal $Sp(2,\mathbb{R})$-bundle. Because $H^1(\mathbb{R}^2,\mathbb{R})=0,$ we know that there exists, up to a bundle isomorphism, only one metaplectic bundle over $\mathbb{R}^2,$ namely the trivial principal $Mp(2,\mathbb{R})$-bundle 
$\mathbb{R}^2 \times Mp(2,\mathbb{R}) \to \mathbb{R}^2$ and thus also a unique metaplectic structure $\Lambda: Mp(2,\mathbb{R}) \times \mathbb{R}^{2} \to Sp(2,\mathbb{R})\times \mathbb{R}^2$ given by
$\Lambda(g,(s,t)):=(\lambda(g),(s,t))$ for $g \in Mp(2,\mathbb{R})$ and $(s,t) \in \mathbb{R}^2.$ Let $\mathcal{S} \to \mathbb{R}^2$ be the symplectic spinor bundle. In this case $\mathcal{S}\to \mathbb{R}^2$ is isomorphic to the trivial vector bundle
${\bf S}\times \mathbb{R}^2 = L^2(\mathbb{R})\times \mathbb{R}^2 \to \mathbb{R}^2.$ Thus, we may think of a symplectic spinor field $\phi$ as of a mapping
$\phi: \mathbb{R}^{2} \to {\bf S} = L^2(\mathbb{R})$. Let us define $\psi: \mathbb{R}^3 \to \mathbb{C}$ by $\psi(s,t,x):=\phi(s,t)(x)$ for each $(s,t,x) \in \mathbb{R}^3.$ One easily shows that $\phi$ is a symplectic Killing spinor if and only if the function $\psi$ satisfies the system

\begin{eqnarray*}
\frac{\partial \psi}{\partial s} &=& \lambda \imath x \psi \mbox{ and}\\
\frac{\partial \psi}{\partial t} &=& \lambda \frac{\partial \psi}{\partial x}.
\end{eqnarray*}

If $\lambda=0,$ the solution of this system of partial differential equations is necessarily
$\psi(s,t,x)=\overline{\psi}(x),$ $(s,t,x) \in \mathbb{R}^3,$ for any $\overline{\psi} \in L^2(\mathbb{R}).$

If $\lambda \neq 0,$ let us consider the independent variable and corresponding dependent variable transformation
$s=s, y=t+\lambda^{-1}x,$ $z=t-\lambda^{-1}x$ and $\psi(s,t,x)=\widetilde{\psi}(s, t+\lambda^{-1}x,t-\lambda^{-1}x)=\widetilde{\psi}(s,y,z).$  
The Jacobian of this transformation is $-2/\lambda \neq 0$ and the transformation is obviously a diffeomorphism.
Substituting this transformation in the studied system, one gets the following equivalent transformed system
\begin{eqnarray*}
\frac{\partial \widetilde{\psi}}{\partial s} &=& \frac{\imath}{2}\lambda^2 (y-z)\widetilde{\psi}\\
\frac{\partial \widetilde{\psi}}{\partial y}+\frac{\partial \widetilde{\psi}}{\partial z}
&=&\lambda(\frac{\partial\widetilde{\psi}}{\partial y}\lambda^{-1}+\frac{\partial \widetilde{\psi}}{\partial z}
(-\lambda^{-1})).
\end{eqnarray*}  
(Let us notice that the substitution we have used is similar to that one which is usually used  to obtain the d'Alemebert's solution of the wave equation in two dimensions.)
The first equation implies $\frac{\partial\widetilde{\psi}}{\partial z}=0,$ and thus
$\widetilde{\psi}(s,y,z)=\overline{\overline{\psi}}(s,y)$ for a function $\overline{\overline{\psi}}.$ Substituting this relation into the second equation of the transformed system, we get
$$\frac{\partial \overline{\overline{\psi}}}{\partial s} = \frac{\imath}{2}(y-z)\lambda^2\overline{\overline{\psi}}.$$
The solution of this equation is
$\overline{\overline{\psi}}(s,y)=e^{\frac{\imath}{2}\lambda^2(y-z)s}\widetilde{\overline{\psi}}(y)$ for a suitable function
$\widetilde{\overline{\psi}}.$ Because of the dependence of the right hand side of the last written equation on $z,$ we see that $\overline{\overline{\psi}}$ does not exist unless 
$\lambda =0$ or $\widetilde{\overline{\psi}}=0$ (More formally, one gets these restrictions by substituting  the last written formula for $\overline{\overline{\psi}}$ into the  first equation of the transformed system.)
Thus, necessarily $\psi =0$ or $\lambda=0.$  The case $\lambda =0$ is excluded by the assumption at the beginning of this calculation.

Summing up, we have proved that any symplectic Killing spinor field $\phi$ on $(\mathbb{R}^2,\omega_0,\nabla)$ is constant, i.e., for each $(s,t) \in \mathbb{R}^2,$ we have $\phi(s,t)=\overline{\psi}$ for a function $\overline{\psi} \in L^2(\mathbb{R}).$  The only symplectic Killing spinor number is zero in this case.

{\bf Remark:} More generally, one can treat the case of a standard symplectic vector space $(\mathbb{R}^{2l}[s^1,\ldots, s^l,t^1,\ldots, t^l],\omega_0)$ equipped with the standard flat Euclidean connection $\nabla.$ One gets by similar lines of reasoning that any symplectic Killing spinor for this Fedosov manifold is also constant, i.e., 
$$\psi(s^1,\ldots, s^l, t^1, \ldots, t^{l}) = \overline{\psi},$$ for 
$(s^1,\ldots, s^l), (t^1,\ldots, t^{l}) \in \mathbb{R}^l$ and  $\overline{\psi} \in L^2(\mathbb{R}^l).$ But we shall see this result more easily\textit{} bellow when we will be studying the prolongated equation mentioned in the Introduction.

Now, in order to make a connection of the symplectic Killing spinor equation to some slightly more known equations, let us introduce the following operators.

The operator
$$\mathfrak{D}:\Gamma(M,\mathcal{S}) \to \Gamma(M,\mathcal{S}), \mbox{   } \mathfrak{D}:=-F^-\nabla^S$$
is called symplectic Dirac operator and its eigenfunctions are called symplectic Dirac spinors.
Let us notice that the symplectic Dirac operator was introduced by Katharina Habermann in 1992. See, e.g., Habermann \cite{KH}.

The operator
$$\mathfrak{T}:\Gamma(M,\mathcal{S})\to \Gamma(M,\mathcal{E}^{11}), \mbox{    } \mathfrak{T} := \nabla^S - p^{10}\nabla^S$$
is called (the first) symplectic twistor operator.

In the next theorem, the symplectic Killing spinor fields are related to the symplectic Dirac spinors and to the kernel of the symplectic twistor operator.

{\bf Theorem  6:} Let $(M,\omega, \nabla)$ be a Fedosov manifold admitting a metaplectic structure.
A symplectic spinor field $\phi \in \Gamma(M,\mathcal{S})$ is a symplectic Killing spinor field if and only if $\phi$ is a symplectic Dirac spinor lying in the kernel of the symplectic twistor operator.

{\it Proof.} We prove this equivalence in two steps.
\begin{itemize}
\item[1)]
Suppose $\phi \in \Gamma(M,\mathcal{S})$ is a symplectic Killing spinor to a symplectic Killing number $\lambda \in \mathbb{C}.$ Thus it satisfies the equation $\nabla^S \phi = \lambda F^+\phi.$
Applying the operator $-F^-$ to the both sides of the preceding equation and using the definition of the symplectic Dirac operator, we get
$\mathfrak{D}\phi = -\lambda F^-F^+\phi = \lambda(-H + F^+ F^-)\phi=-\lambda H \phi = - \lambda (-\imath l \phi) = \imath \lambda l \phi$  due to the definition of $H$ and the Lemma 3 (part 2a and 3). Thus $\phi$ is a symplectic Dirac spinor. 

Now, we compute $\mathfrak{T}\phi.$ Using the definition of $\mathfrak{T}$, we get $\mathfrak{T}\phi = (\nabla^S-p^{10}\nabla^S)\phi=
\lambda (F^+\phi - p^{10}F^+\phi) = \lambda p^{11} F^+ \phi = 0,$ because $F^+ \phi \in \Gamma(M,\mathcal{E}^{10})$ due to the Lemma 3 part 2a.

\item[2)]
Conversely, let $\phi \in \Gamma(M,\mathcal{E}^{00})$ be in the kernel of the symplectic twistor operator and also a symplectic Dirac spinor. Thus, we have $\nabla^S \phi - p^{10} \nabla^S \phi =0$ and $\mathfrak{D}\phi=-F^-\nabla^S \phi = \mu \phi$ for a complex number $\mu \in \mathbb{C}.$ 
From the first  equation, we deduce that $\psi := \nabla^S \phi \in \Gamma(M,\mathcal{E}^{10}).$
Because $F^+_{|\Gamma(M,\mathcal{E}^{00})}$ is surjective onto $\Gamma(M,\mathcal{E}^{10})$ (see Lemma 3 part 2b),
there exists a $\psi'\in \Gamma(M,\mathcal{E}^{00})$ such that $\psi=F^+\psi'.$ Let us compute $F^+F^-\psi = F^+F^-F^+\psi'=F^+(H - F^+F^-)\psi'=F^+(-\imath l \psi') =-\imath l \psi,$ where we have used the defining equation for $H$ and the Lemma 3 part 2a and the part 3. Thus we get
\begin{eqnarray}
-F^+F^-\psi=\imath l \psi.
\end{eqnarray}
 From the  symplectic Dirac equation, we get $\mu \phi = -F^-\psi.$ Thus $-F^+F^-\psi=\mu F^+\phi.$ Using the equation (1), we obtain
$\imath l \psi = \mu F^+ \phi,$ i.e., $\nabla^S\phi =  -\imath\frac{\mu}{l}F^+\phi.$
Thus, $\phi$ is a symplectic Killing spinor to the symplectic Killing spinor number $-\imath \mu/l.$
\end{itemize} $\Box$

In the next theorem, we derive the mentioned prolongation of the symplectic Killing spinor equation.
It is a zeroth order equation. More precisely, it is an equation for the sections of the kernel of an endomorphism of the symplectic spinor bundle $\mathcal{S}\to M.$ A similar computation  is well known from the Riemannian spin geometry. See, e.g., Friedrich \cite{Friedrich}.

{\bf Theorem 7:} Let $(M^{2l},\omega,\nabla)$ be a Fedosov manifold admitting  a metaplectic structure and a symplectic Killing spinor field $\phi \in \Gamma(M,\mathcal{S})$ to the symplectic Killing spinor number $\lambda$.
Then 

      $$\sigma^{ij}e_{i}.e_{j}.\phi = 2 l \lambda^2 \phi.$$ 

{\it Proof.} Let $\phi \in \Gamma(M^{2l},\mathcal{S})$ be a symplectic spinor Killing field, i.e., $\nabla^S_X\phi=\lambda X. \phi$ for a complex number $\lambda$ and any vector field $X\in \mathfrak{X}(M).$ For vector fields $X, Y \in \mathfrak{X}(M),$ we may write 
\begin{eqnarray*} 
R^S(X,Y)\phi&=&(\nabla_X\nabla_Y - \nabla_Y\nabla_X - \nabla_{[X,Y]})\phi\\
            &=&\lambda \nabla_X (Y.\phi) - \lambda \nabla_Y (X.\phi) - \lambda[X,Y].\phi\\
            &=&\lambda (\nabla_XY).\phi + \lambda Y. (\nabla_X\phi) - \lambda (\nabla_YX). \phi - \lambda X. \nabla_Y.\phi - \lambda [X,Y].\phi\\
            &=&\lambda T(X,Y).\phi + \lambda^2 (Y.X.-Y.X.)\phi\\
            &=&\lambda T(X,Y) .\phi  + \imath \lambda^2 \omega(X,Y)\phi = \imath \lambda^2 \omega(X,Y)\phi,
\end{eqnarray*}            
where we have used the symplectic Killing spinor equation and the compatibility of the symplectic spinor covariant derivative and the symplectic Clifford multiplication (Lemma 4).          
     
Thus $R^S \phi = \imath \lambda^2 \omega \otimes \phi.$ Because of the Lemma 3 part 2c,
we know that the right hand side is in $\Gamma(M,\mathcal{E}^{20}).$ Thus also $R^S\phi = p^{20} R^S \phi.$ 
Using the Lemma 5,
we get $\frac{\imath}{2l}\omega \otimes \sigma^{ij}e_{i}.e_{j}.\phi = \imath \lambda^2 \omega \otimes \phi.$ 
Thus $\sigma^{ij}e_{i}.e_{j}.\phi = 2 l \lambda^2 \phi$ and the theorem follows.  
$\Box$

{\bf Remark:} Let us recall that in the Riemannian spin geometry (positive definite case), the existence of a non-zero Killing spinor implies that the manifold is Einstein. 
Further, let us notice that if the symplectic Ricci curvature tensor $\sigma$ is (globally) diagonalizable by a symplectomorphism, the prolongated equation has the shape of the equation for eigenvalues of the Hamiltonian of an elliptic $l$ dimensional harmonic oscillator with possibly varying axes lengths. An  example of a diagonalizable symplectic Ricci curvature will be treated in the Example 3.  Although, in this case the axis will be constant and the harmonic oscillator will be spherical.




Now, we derive a simple consequence of the preceding theorem in the case of Fedosov manifolds of Weyl type, i.e., $\sigma=0$.

{\bf Corollary 8:} Let $(M,\omega,\nabla)$ be a Fedosov manifold of Weyl type.
Let $(M,\omega)$ admits a metaplectic structure and a symplectic Killing spinor $\phi$ field to the symplectic Killing spinor number $\lambda.$
Then the symplectic Killing spinor number $\lambda=0$ and $\phi$ is locally covariantly  constant.

{\it Proof.} Follows immediately from the preceding theorem and the symplectic Killing spinor equation. $\Box$

{\bf Example 2:} Let us go back to the case of  $(\mathbb{R}^{2l},\omega_0,\nabla)$ from the Remark bellow the Example 1.
The Corrollary 8 implies that any symplectic Killing spinor field for this structure is covariantly constant, i.e., in fact constant in this case, and any symplectic Killing number is zero. In this case, we see that the prolongated equation from the Theorem 7 makes it able to compute the symplectic Killing spinor fields without any big effort, compared to the calculations in the Example 1 where the $2$-plane was treated.

In the next example, we compute the symplectic Killing spinor fields on $S^2$ equipped with the standard symplectic structure and the Riemannian connection of the round metric. This is an example of a Fedosov manifold (specified more carefully bellow) for which one can not use the Corollary 8, because it is not of Weyl type. But still, one can use the Theorem 7.

{\bf Example 3:} Consider the round sphere $(S^2,r^2(d\theta^2+\sin^2 \theta d\phi^2))$ of radius $r>0,$ $\theta$ being the longitude a $\phi$ the latitude. Then $\omega:=r^2\sin \theta d\theta \wedge d\phi$ is the volume form of the round sphere. Because $\omega$ is also a symplectic form, $(S^2,\omega)$ is a symplectic manifold. Let us consider the Riemannian  connection $\nabla$ of the round sphere. Then $\nabla$ preserves the symplectic volume form $\omega$ being a metric connection of the round sphere. Because $\nabla$ is torsion-free, we see that $(S^2,\omega, \nabla)$ is a Fedosov manifold. 
Now, we will work in a coordinate patch without mentioning it explicitly.
Let us set $e_1:=\frac{1}{r}\frac{\partial}{\partial \theta}$ and $e_2:=\frac{1}{r \sin \theta}\frac{\partial}{\partial \phi}.$ 
Clearly, $\{e_1,e_2\}$ is a local adapted symplectic frame and it is a local orthogonal frame as well.
With respect to this basis,
the Ricci form $\sigma$ of $\nabla$ takes the form

$$[\sigma^{ij}]_{i,j=1,2}
= 
\left( \begin{array}{cc}
1/r & 0\\
0 & 1/r \end{array} 
\right).$$

Let us consider $S^2$ as the complex projective space $\mathbb{CP}^1.$ It is easy to see that the (unique) complex structure on $\mathbb{CP}^1$ is compatible with the volume form. The first Chern class of the tangent bundle to $\mathbb{CP}^1$ is known to be even. Thus, the symplectic manifold $(S^2,\omega)$ admits a metaplectic structure. Thus we may consider a symplectic Killing spinor field $\phi \in \Gamma(S^2, \mathcal{S})$ corresponding to a symplectic Killing spinor number $\lambda.$
Because the first homology group of the sphere $S^2$ is zero, the metaplectic structure is unique and thus the trivial one.
 Because of the triviality of the associated symplectic spinor bundle $\mathcal{S} \to S^2,$ we may write $\phi(m)=(m,f(m))$ where $f(m) \in L^2(\mathbb{R})$ for each $m \in S^2.$
Using the Theorem 7 and the prescription for the Ricci form, we get that
$\sigma^{ij}e_i.e_j.[f(m)]=\frac{1}{r}H[f(m)]=2\lambda^2 f(m),$ where
$H = \frac{\partial^2}{\partial x^2} - x^2$ is the quantum Hamiltonian of the one dimensional harmonic oscillator.
The solutions of the Sturm-Liouville type equation $H[f(m)] = 2r\lambda^2 f(m),$ $m \in S^2,$ are well known.
The eigenfunctions of $H$ are the Hermite functions
$f_l(m)(x)=h_l(x):=e^{x^2/2}\frac{d^l}{dx^l}(e^{-x^2})$ for $m \in S^2$ and $x\in \mathbb{R}$ and the corresponding
eigenvalues are $-(2l+1),$ $l \in \mathbb{N}_0.$
Thus $2 r \lambda^2 = -(2l+1)$ and consequently
$$\lambda = \pm \imath \sqrt{\frac{2l+1}{2r}}.$$

Using the fact that $\{e_1,e_2\}$ is a local orthonormal frame and $\nabla$ is metric and torsion-free, we easily get
$$\begin{array}{ll}
\nabla_{e_1}e_1=0 & \nabla_{e_1}e_2=0\\
\nabla_{e_2}e_1=\frac{\cot \theta}{r}e_2&  \nabla_{e_2}{e_2}=-\frac{\cot \theta}{r}e_1.
\end{array}$$

From the definition of differentiability of functions with values in a Hilbert space, 
we see easily as a consequence of the preceding computations that any symplectic Killing spinor field is necessarily of the form
$\phi(m)=(m, c(m)f_l(m))$ for a smooth function $c\in \mathcal{C}^{\infty}(S^2,\mathbb{C}).$
Substituting this Ansatz into the symplectic Killing spinor equation,
 we get for each vector field $X\in \mathfrak{X}(S^2)$ the equation
\begin{eqnarray*}
\nabla_X(c f_l) = (Xc)f_l+c\nabla_X f_l = \lambda c(X.f_l).
\end{eqnarray*}

Due to the Lemma 4, we have for a local adapted symplectic frame $s:U\subseteq S^2 \to \mathcal{P}=Sp(2,\mathbb{R})\times S^2,$
$$\nabla_X f_l = [\overline{s}, X(f_l)_s] - \frac{\imath}{2}[e_2.(\nabla_X e_1). -e_1.(\nabla_X e_2).]f_l$$
(See the paragraph above the Lemma 4 for an explanation of the notation used in this formula.)

Because $m\mapsto (m,f_l(m))$ is constant as a section of the trivial bundle
$\mathcal{S} \to S^2,$ the first summand of the preceding expression vanishes. Thus for $X =e_1,$ we get
$$(e_1 c)f_l + \frac{\imath c}{2}[e_2.(\nabla_{e_1}e_1). - e_1.(\nabla_{e_1} e_2).]f_l = \lambda c (e_1.f_l).$$
Using the knowledge of the values of $\nabla_{e_1}e_j,$ for $j=1,2,$ computed above, the second summand at the left hand side of the last written equation vanishes and thus, we get
$$\frac{1}{r}\frac{\partial c}{\partial \theta} f_l = \lambda c \imath x f_l.$$

This equation implies
$c(\theta, \phi)= \psi(x,\phi)e^{\imath r x \lambda \theta}$
for $x$ such that $h_l(x)\neq 0$ and a suitable function $\psi.$ (The set of such $x \in \mathbb{R},$ such that $h_l(x)\neq 0$ is the complement in $\mathbb{R}$ of a finite set.)
Because $r>0$ is given and $\lambda$ is certainly non-zero (see the prescription for $\lambda$ above), the only possibility for $c$ to be independent of $x$ is
$\psi=0.$ Therefore $c=0$ and consequently $\phi=0.$  On the other hand, $\phi=0$ (the zero section) is clearly a solution, but according to the definition not a symplectic Killing spinor. Thus, there is no symplectic Killing spinor field on the round sphere.

{\bf Remark:} In the future, one can study holonomy restrictions implied by the existence of a symplectic Killing spinor. One can also try to extend the results to general symplectic connections, i.e., to drop the condition on the torsion-freeness or  study also the symplectic Killing fields on Ricci type Fedosov manifolds admitting a metaplectic structure in more detail.

\end{document}